\documentclass[12pt,czech]{amsart}
\usepackage{amsmath}
\usepackage{amscd}
\usepackage[psamsfonts]{amssymb}
\usepackage[mathscr,psamsfonts]{eucal}

\DeclareMathAlphabet{\eurm}{U}{eur}{m}{n}
\newcommand{\tfr}[2]{\tfrac{#1}{#2}\,}
\newcommand{\F}{{\mathcal F}(Y_1,Y_2)}
\newcommand{\Sec}{{S}(E,Q)}
\newcommand{\mc}[1]{{\mathcal{#1}}}
\newcommand{\ucar}[1]{\underset{#1}{\times}}
\newcommand{\der}{\partial}
\newcommand{\wed}{\wedge}

\DeclareMathOperator{\id}{{id}}
\DeclareMathOperator{\Aut}{{Aut}}
\DeclareMathOperator{\byd}{\,{\raisebox{.1ex}{$\eurm :$}{\eurm =}}\,}

\newcommand{\sig}{\sigma}
\newcommand{\alp}{\alpha}
\newcommand{\lam}{\lambda}
\newcommand{\bet}{\beta}
\newcommand{\del}{\delta}
\newcommand{\gam}{\gamma}
\newcommand{\bEq}{\begin{eqnarray}}
\newcommand{\eEq}{\end{eqnarray}}
\newcommand{\beq}{\begin{eqnarray*}}
\newcommand{\eeq}{\end{eqnarray*}}

\newcommand{\com}{\circ}
\newcommand{\car}{\times}
\newcommand{\ten}{\otimes}
\newcommand{\Rn}{{\mathbb R}}
\newcommand{\mto}{\mapsto}

\begin{document}

\font\aa=msam10 scaled \magstep 4

\newcounter{theorem}

\newtheorem{definition}[theorem]{Definition}
\newtheorem{lemma}[theorem]{Lemma}
\newtheorem{proposition}[theorem]{Proposition}
\newtheorem{theorem}[theorem]{Theorem}
\newtheorem{corollary}[theorem]{Corollary}
\newtheorem{remark}[theorem]{Remark}
\newtheorem{example}[theorem]{Example}
\newtheorem{Note}[theorem]{Note}
\newcounter{assump}
\newtheorem{Assumption}{\indent Assumption}[assump]
\renewcommand{\thetheorem}{\thesection.\arabic{theorem}}

\newcommand{\bCr}{\begin{corollary}}
\newcommand{\eCr}{\end{corollary}}
\newcommand{\bDf}{\begin{definition}\em}
\newcommand{\eDf}{\end{definition}}
\newcommand{\bLm}{\begin{lemma}}
\newcommand{\eLm}{\end{lemma}}
\newcommand{\bPr}{\begin{proposition}}
\newcommand{\ePr}{\end{proposition}}
\newcommand{\bRm}{\begin{remark}\em}
\newcommand{\eRm}{\end{remark}}
\newcommand{\bEx}{\begin{example}\em}
\newcommand{\eEx}{\end{example}}
\newcommand{\bTh}{\begin{theorem}}
\newcommand{\eTh}{\end{theorem}}
\newcommand{\bNt}{\begin{Note}\em}
\newcommand{\eNt}{\end{Note}}
\newcommand{\bPf}{\begin{proof}[\noindent\indent{\sc Proof}]}
\newcommand{\ePf}{\end{proof}}
\title[On the geometry of variational calculus]
        {On the geometry of variational calculus\\ on some
        functional bundles}

\author[A. Cabras, J. Jany\v ska, I. Kol\'a\v{r}]
        {Antonella Cabras, Josef Jany\v ska, Ivan Kol\'a\v{r}
}




\address{
\newline
Antonella Cabras
\newline
Department of Applied Mathematics, Florence University
\newline
Via S. Marta 3, 50139 Florence, Italy
\newline
E-mail: {\tt cabras@dma.unifi.it}
\newline
{\ }
\newline
Josef Jany\v{s}ka
\newline
Department of Mathematics, Masaryk University
\newline
Jan\'a\v ckovo n\'am. 2a, 662 95 Brno, Czech Republic
\newline
E-mail: {\tt janyska@math.muni.cz}
\newline
{\ }
\newline
Ivan Kol\'a\v{r}
\newline
Department of Algebra and Geometry, Masaryk University
\newline
Jan\'a\v ckovo n\'am. 2a, 662 95 Brno, Czech Republic
\newline
E-mail: {\tt kolar@math.muni.cz}
}

\thanks{The second author was supported by the
Grant Agency of the Czech Republic
under the project number GA 201/02/0225 and
the third author was supported
by the Ministry of Education of the Czech Republic under
the project
MSM 143100009.}

\begin{abstract}
We first generalize the operation of formal exterior
differential in the case of finite dimensional fibered manifolds
and then we extend it to certain bundles of smooth maps. In order to
characterize the operator order of some morphisms between our bundles
of smooth maps, we introduce the concept of fiberwise $(k,r)$-jet.
The relations to the Euler-Lagrange morphism of the variational
calculus are described.
\end{abstract}
\maketitle

\section*{Introduction}
\setcounter{equation}{0}
\setcounter{theorem}{0}

Our geometrical research was inspired by the paper on the
Schr\"odinger operator by the second author and M. Modugno,
\cite{JanMod02}, as well as by their previous joint paper with
A. Jadczyk, \cite{JadJanMod98}. We are interested mainly in certain
geometric objects and operations related with the functional bundle
${S}(E,Q)$ of all sections $E_x\to Q_x$ of a 2-fibered manifold
$Q\to E\to M$, $x\in M$.
In \cite{CabKol95}, the first and the third authors established
the theory of connections in a somewhat more general situation
of the bundle $\F\to M$ of all smooth maps between the fibers
over the same base point of two fibered manifolds
$Y_1\to M$ and $Y_2\to M$ with the same base $M$. The main purpose
of the present paper is to introduce some geometric concepts and
to study some geometric operations that could be useful
for the variational calculus on these functional bundles.

Our approach to the variational calculus is based on the
formal exterior differential on finite dimensional fibered
manifolds introduced by A. Trautman, \cite{Tra72}, and further
developed by the third author, \cite{Kol84, Kol84a}.
In Section 1 of the present paper we study a slight finite
dimensional generalization of this concept in a form
suitable for our next purposes. Section 2 is devoted
to some geometric properties of the bundles $\F$ and $\Sec$
in the framework of the Fr\"olicher's theory of smooth
structures, \cite{Fro82}. The morphisms between our functional
bundles represent a kind of differential operators.
As pointed out already in \cite{CabKol95}, one can distinguish
an important class of them that have finite order
in the operator sense. In Section 3 we modify this idea
to the morphisms defined on the $r$-th jet prolongation $J^r\F$.
This leads us to an original concept of fiberwise $(k,r)$-jet
of a base preserving morphism of finite dimensional
fibered manifolds.
Section 4 deals with the formal exterior differentiation
over the functional bundle $\F$. In Section 5 we study its
restriction to the bundle $\Sec$ of sections.
In Proposition \ref{Pr3}
we characterize an important situation in which the finite
dimensional formal exterior differential and the analogous
operation over $\Sec$ are naturally related.
Finally, Section 6 is devoted to the Euler-Lagrange morphism
on $\Sec$ from the viewpoint of our previous operations.

If we deal with finite dimensional manifolds and maps
between them, we always assume they are of class $C^\infty$, i.e.
smooth in the classical sense. On the other hand, the smooth
spaces and maps in the sense of A. Fr\"olicher are said to be
$F$-smooth. Unless otherwise specified, all morphisms
are assumed to be base preserving. In all standard situations
we use the terminology and notation from the monograph
\cite{KolMicSlo93}.

We acknowledge Marco Modugno for suggesting this subject
and for several stimulating discussions.

\section{The formal exterior differential in finite dimension}
\setcounter{equation}{0}
\setcounter{theorem}{0}

We recall that $Q\overset{q}{\to} E \overset{p}{\to} M$ is said to be a
2-fibered manifold, if both $q$ and $p$ are surjective submersions.
Consider two 2-fibered manifolds
$Z\to Y\to M$, $W\to Y\to M$, a fibered manifold $N\to M$ and
a morphism
\beq
    \psi:Z \ucar{Y} W\to N\,,
\eeq
where $Z\ucar{Y} W$ is interpreted as a fibered manifold over $M$.
Then the rule
\bEq\label{Eq1}
    J^k\psi(j^k_x s, j^k_x\sigma)= j^k_x\psi(s,\sigma)
\eEq
defines a map
\bEq\label{Eq2}
    J^k \psi:J^kZ\ucar{J^k Y} J^k W\to J^k N\,.
\eEq
In the case of
\beq
    \psi: J^r Z\ucar{J^t Y} J^s W\to  N\,, \quad r\ge t\le s\,,
\eeq
we obtain
\beq
    J^k \psi:J^kJ^r Z\ucar{J^kJ^t Y} J^kJ^s W\to J^k N\,.
\eeq
Then we introduce
\bEq\label{Eq3}
    J^k_{\mathrm{hol}}\psi:J^{k+r} Z\ucar{J^{k+t} Y} J^{k+s} W
        \to J^k N
\eEq
by means of the canonical inclusions of the holonomic jet prolongations
into the iterated jet prolongations.

In particular, consider
\bEq\label{Eq4}
    \varphi:J^r Y\ucar{J^s Y} VJ^s Y\to Z\,,\quad s\le r\,,
\eEq
where $Z\to M$ is a fibered manifold. By using the well known
identification $\varkappa_s: VJ^s Y\to J^sVY$, we construct
\beq
\varphi\com(\id_{J^r Y}\ucar{J^s Y}\varkappa^{-1}_s):
        J^rY\ucar{J^s Y} J^sVY\to Z
\eeq
and
\beq
J^k_{\mathrm{hol}}(\varphi\com(\id_{J^r Y}\ucar{J^s Y}
                \varkappa^{-1}_s)): J^{k+r}Y\ucar{J^{k+s} Y}
        J^{k+s}VY\to J^k Z\,.
\eeq

\noindent
Then we define
\begin{align}\label{Eq5}
    \mc{J}^k_{\mathrm{hol}}\varphi & := J^k_{\mathrm{hol}}
        (\varphi\com(\id_{J^r Y}\ucar{J^s Y}\varkappa^{-1}_s))
        \com(\id_{J^{k+r} Y}\ucar{J^{k+s} Y}\varkappa_{s+k}):
\\
  &   :   J^{k+r} Y \ucar{J^{k+s} Y} VJ^{k+s} Y\to J^k Z\,.\nonumber
\end{align}

Let $\eta:Y\to VY$ be a vertical vector field on $Y$ and
$\mc{J}^s\eta:J^s Y\to VJ^s Y$ be its flow prolongation. Write
\beq
    \varphi(\mc{J}^s\eta)=\varphi\com ( \id_{J^r Y}\ucar{J^s Y}
        \mc{J}^s \eta): J^r Y\to Z\,.
\eeq
Then
\beq
    J^k_{\mathrm{hol}}(\varphi(\mc{J}^s\eta)):J^{k+r} Y\to J^k Z\,.
\eeq
On the other hand,
\beq
  \mc{J}^k_{\mathrm{hol}}\varphi:
        J^{k+r} Y \ucar{J^{k+s} Y} VJ^{k+s} Y\to J^k Z\,,
\eeq
so that
\beq
        (\mc{J}^k_{\mathrm{hol}}\varphi)(\mc{J}^{k+s}\eta )
: J^{k+r} Y\to J^k Z\,.
\eeq

\bPr\label{Pr1}
For every $\varphi$ and $\eta$, we have
\bEq\label{Eq6}
        (\mc{J}^k_{\mathrm{hol}}\varphi)(\mc{J}^{k+s}\eta )=
        J^{k}_{\mathrm{hol}}(\varphi(\mc{J}^s\eta))\,.
\eEq
\ePr

\bPf
This follows from the well known fact
$\mc{J}^s\eta = \varkappa^{-1}_s\com J^s\eta$, where
$J^s\eta:J^s Y\to J^s VY$ is the functorial prolongation of $\eta$.
\ePf

Consider the case $Z=\bigwedge^l T^* M$ in (\ref{Eq4}).
The exterior differential $d$ on $M$ is a first order operator,
so that $d$ determines the associated map
$\del:J^1\bigwedge^lT^*M\to \bigwedge^{l+1}T^*M$ satisfying
$d\omega = \del\com (J^1\omega)$ for every $l$-form
$\omega:M\to \bigwedge^lT^*M$.

\bDf\label{Df1}
For every morphism
$\varphi:J^r Y\ucar{J^{s} Y} VJ^{s}Y\to \bigwedge^l T^*M$,
we define its {\em formal exterior differential} by
\bEq\label{Eq7}
D\varphi := \del\com (\mc{J}^1_{\mathrm{hol}}\varphi):
J^{r+1} Y\ucar{J^{s+1} Y} VJ^{s+1} Y\to \bigwedge^{l+1}T^* M\,.
\eEq
\eDf

Proposition \ref{Pr1} implies that this concept represents
a generalization of that one introduced by the third author in
\cite{Kol84}. In fact, $\varphi$ is assumed to be linear
in $VJ^{s} Y$ in \cite{Kol84},
while in (\ref{Eq7}) $\varphi$ is quite arbitrary.

Consider some local fiber coordinates $x^i,x^p$ on $Y$,
$i=1,\dots,m=\dim M$, $p=m+1,\dots, m+n = \dim Y$.
Let $\alp$ and $\sig$ be multiindices of the range $m$. Write
\beq
    x^p_\alp\,,\quad 0\le \parallel\!\! \alp \!\! \parallel \le r
\eeq
for the induced coordinates on $J^r Y$ and
\beq
        x^p_\sig\,,\,\, X^p_\sig= dx^p_\sig, \quad 0\le ||\sig|| \le s
\eeq
for the induced coordinates on $VJ^sY$. If
\bEq\label{Eq1.8}
    a_{i_1\dots i_l}(x^i,x^p_\alp, X^p_\sig)\, dx^{i_1}\wed \dots
        \wed dx^{i_l}
\eEq
is the coordinate expression of $\varphi$, then the coordinate
form of $D\varphi$ is
\bEq\label{Eq1.9}
    \left( \frac{\der a_{i_1\dots i_l}}{\der x^i } +
        \frac{\der a_{i_1\dots i_l}}{\der x^p_\alp }\, x^p_{\alp i} +
        \frac{\der a_{i_1\dots i_l}}{\der X^p_\sig }\, X^p_{\sig i}\right)
        \, dx^i\wed dx^{i_1}\wed \dots \wed dx^{i_l}\,.
\eEq

\section{The functional bundle $\Sec$}
\setcounter{equation}{0}
\setcounter{theorem}{0}

We shall use the following simplified version, \cite{CabJanKol04},
of the theory of smooth spaces by A. Fr\"olicher, \cite{Fro82}.
An $F$-smooth space is a set $S$ along with a set $C_{S}$
of maps $\gam:\Rn\to S$, which are called $F$-smooth curves,
satisfying

(i) each constant curve $\Rn \to S$ belongs to $C_{S}$,

(ii) if $\gam\in C_{S}$ and $\varepsilon\in C^\infty(\Rn,\Rn)$,
then $\gam\com \varepsilon\in C_{S}$.

\noindent
Every subset $\bar S\subset S$ is also an $F$-smooth space,
if we define $C_{\bar S}\subset C_S$ to be the subset of all curves
with values in $\bar S$.
If $(S',C_{S'})$ is another $F$-smooth space, a map
$f:S\to S'$ is said to be $F$-smooth, if $f\com \gam$ is an
$F$-smooth curve on $S'$ for every $F$-smooth curve $\gam$ on $S$.
So we obtain the category $\mc{S} $ of $F$-smooth spaces.

In particular, every smooth
manifold $M$ turns out to be an $F$-smooth space by assuming
as $F$-smooth curves just the smooth curves. Moreover, a map between
smooth manifolds is $F$-smooth, if and only if it is smooth.
An $F$-smooth bundle is a triple of an $F$-smooth space $S$,
a smooth manifold $M$ and a surjective $F$-smooth map
$p:S\to M$. If $p':S'\to M'$ is another
$F$-smooth bundle, then a morphism of $S$ into $S'$ is a pair of
an $F$-smooth map $f:S\to S'$ and a smooth map $\underline f:
M\to M'$ satisfying $\underline f\com p= p'\com f$.
So we obtain the category ${\mc{S}\mc{B}}$ of
$F$-smooth bundles.

If $p_1: Y_1\to M$,  $p_2: Y_2\to M$ are two fibered manifolds,
we write
\beq
\F=\underset{x\in M}{\bigcup}
        C^\infty (Y_{1x},Y_{2x})
\eeq
and denote by $p:\F\to M$ the canonical
projection. A curve $\widehat c:\Rn\to \F$ is said to be
$F$-smooth, if $\underline{c}\byd p\com\widehat c:\Rn\to M$ is a smooth
curve and the induced map
\beq
{c}:\underline{c}^* Y_1\to Y_2,\quad {c}(t,y)=
        \widehat c(t)(y),\quad p_1(y)=\underline{c}(t)\,,
\eeq
is also smooth, \cite{CabKol95}. The $F$-smooth sections of $\F$
are identified with the base preserving morphisms $s:Y_1\to Y_2$.
We write $\widehat s:M\to \F$ for the $F$-smooth section induced by $s$.

The tangent bundle $T\F\to TM$ is defined as follows,
\cite{CabKol95}. For every $F$-smooth curve $\widehat{f}:\Rn \to \F$,
we first construct the tangent vector
$X=\frac{\der}{\der t}|_0(p\com \widehat{f})\in TM$. Write
\beq
    T_XY_1= (Tp_1)^{-1} (X)\subset TY_1, \quad T_X Y_2
        = (Tp_2)^{-1}(X)\subset TY_2\,.
\eeq
Then $\widehat f$ defines a map $T_0\widehat f:T_X Y_1\to T_X Y_2$ by
\bEq\label{Eq8}
    T_0\widehat f(\frac{\der}{\der t}|_0 h(t))
        = \frac{\der}{\der t}|_0\widehat f(t)(h(t))\,,
\eEq
where we may assume that $h:\Rn \to Y_1$ satisfies
$p\com \widehat f= p_1\com h$. We say that $\widehat f$ and another
$F$-smooth curve $\widehat g:\Rn\to \F$ satisfying
$\frac{\der}{\der t}|_0 (p\com \widehat g) = X$ determine
the same tangent vector at
$f(0)=g(0)\in \F$, if $T_0\widehat f= T_0 \widehat g:T_XY_1\to T_X Y_2$.
The set $T\F$ of all equivalence classes is called the
tangent bundle of $\F$. The map $T_0\widehat f$ is said to be
the associated map of the tangent vector $\frac{d}{dt}|_0\widehat f$.

Since $T\F\subset \mc{F}(TY_1\to TM, TY_2\to TM)$, this is also an
$F$-smooth bundle. The vertical tangent bundle $V\F\to M$
is the subbundle of $T\F$ of all elements projected by $Tp$ into
a zero vector on $M$.

Given a 2-fibered manifold
$Q\overset{q}{\to} E\overset{p}{\to} M$,
we denote by $\Sec\subset\mc{F}(E,Q\to M)$ the $F$-smooth bundle
of all sections $s:E_x\to Q_x$ of $q$.

A 2-fibered manifold morphism is a triple $(f,f_1,f_0)$ such that
the following diagram commutes
$$
\begin{CD}
        Q @>q>> E @>p>> M
\\
        @VfVV @Vf_1VV @Vf_0VV
\\
        \bar Q @>\bar q>> \bar E @>\bar p>> \bar M
\end{CD}
$$
So we obtain the category $2\mc{FM}$. Write $2\mc{FM}^I\subset
2\mc{FM}$ for the category defined by the requirement that $f_1$ is a
diffeomorphism on each fiber. If $f\in 2\mc{FM}^I$, we have
the induced map
\beq
    S(f):S(E,Q)\to S(\bar E,\bar Q)
\eeq
transforming $s:E_x\to Q_x$ into
\beq
    f_x\com s \com (f_{1x})^{-1}:\bar E_{f_0(x)}\to \bar Q_{f_0(x)}\,.
\eeq
Clearly, $S$ is a functor on $2\mc{FM}^I$ with values in $\mc{SB}$.

If we have another 2-fibered manifold $P\to E\to M$,
then a $2\mc{FM}$-morphism over $\id_E$ will be called
an $E$-morphism. In this case we shall also write
$\widehat f= S(f):S(E,Q)\to S(E,P)$.

Consider a vertical curve $\widehat\gam:\Rn\to \Sec$
over $x\in M$. Then $\gam(t):E_x\to Q_x$, $t\in\Rn$,
and $\gam(t)(y)$ is a vertical
curve on $Q\to E$ for every $y\in E_x$.
Hence $\frac{d}{dt}|_0\gam(t)(y)\in V_y(Q\to E)$.
Using the standard globalization procedure, \cite{Slo86}, we deduce
\bEq\label{Eq9}
    V\Sec={S}(E,V(Q\to E))\,.
\eEq

We have a canonical injection
\bEq\label{Eq10}
    i:J^r\Sec\to {S}(E,J^r(Q\to E))
\eEq
defined as follows. Consider a section $\widehat s:M\to \Sec$,
so that $s:E\to Q$. Then $j^r\widehat s$ determines $j^rs:E\to J^r(Q\to E)$.
We have $j^r_x \widehat s\in J^r_x\Sec\subset J^r_x\mc{F}(E,Q)$ and we set
\beq
    i(j^r_x \widehat s) = j^rs|E_x:E_x\to J^r_x(Q\to E)\,.
\eeq
We shall consider some local fiber coordinates $x^i, x^p$ on
$Y_1$ and   local fiber coordinates $x^i,z^a$
on $Y_2$.
In the case of $Q\to E\to M$, $x^i,x^p,z^a$ will mean the corresponding
fiber coordinates on $Q$. Hence the coordinate expression of
$j^r_x\widehat s$ are the functions
\beq
    z^a_\alp(x^p),\quad 0\le \parallel\!\!\alp\!\!\parallel\le r\,,
\eeq
where $\alp$ is a multiindex of the range $m$, \cite{CabKol95}.
On the other hand, the coordinate expression of $i(j^r_x \widehat s)$ are
some functions $z^a_{\alp\bet}(x^p)$,
$0\le \parallel\!\!\alp\!\!\parallel+\parallel\!\!\bet\!\!\parallel\le r$,
where $\bet$ is a multiindex of the range $m+1,\dots, m+n$.
Our definition implies
\bEq\label{Eq11}
    z^a_{\alp\bet}=\der_\bet z^a_\alp(x^p)\,.
\eEq

\bRm\label{Rm2.1}
We remark that (\ref{Eq11}) describes also a general injection
\bEq\label{Eq12}
    {S}(E,J^r(Q\to B))\hookrightarrow {S}(E,J^r(Q\to E))\,.
\eEq
\eRm

Consider another 2-fibered manifold $P\to E\to M$ and
an $E$-morphism $f:Q\to P$. Then we have the induced maps
\beq
    J^rf :J^r(Q\to E)\to J^r(P\to E),\quad J^r f(j^r_y s)=j^r_y(f\com s)
\eeq
and ${S}(f):\Sec\to {S}(E,P)$.

\bLm\label{Lm1}
The following diagram commutes
\bEq\label{Eq13}
\begin{CD}
    J^r\Sec @>J^r {S}(f)>> J^r {S}(E,P)
\\
        @ViVV @ViVV
\\
{S}(E,J^r(Q\to E)) @>{S}(J^rf)>> {S}(E,J^r(P\to E))
\end{CD}
\eEq
\eLm

\bPf
For $j^r_x \widehat s\in J^r \Sec$, we obtain clockwise
$i(j^r_x(\widehat f\com \widehat s))=i(j^r_x(\widehat{f\com s}))
= j^r(f\com s)|E_x$
and counterclockwise ${S}(J^rf)(j^rs|E_x)=j^r(f\com s)|E_x$.
\ePf

The classical exchange map $\varkappa_r:VJ^rY\to J^rVY$ is
defined by
\beq
    \tfr{\der}{\der t}|_0j^r_x s(t,u)
\mto j^r_x \tfr{\der}{\der t}|_0s(t,u),\quad t\in \Rn, u\in M\,,
\eeq
\cite{KolMicSlo93}. In the functional case, we have an exchange map
\beq
    K_r:VJ^r\F\to J^rV\F
\eeq
defined by the analogous formula
\bEq\label{Eq14}
    K_r\big(\tfr{\der}{\der t}|_0j^r_x\widehat s(t,u)\big) =
        j^r_x\tfr{\der}{\der t}|_0\widehat s(t,u)\,.
\eEq
If we consider $\Sec$ instead of $\F$, then the values of $\widehat s$
in (\ref{Eq14}) are the sections of $q$, so that
we have a restricted and corestricted map
\beq
    K_r:VJ^r\Sec\to J^r V\Sec\,.
\eeq
The same character of the definitions of $\varkappa_r$ and $K_r$ implies
that the following diagram commutes
\bEq\label{Eq15}
\begin{CD}
    VJ^r\Sec @>K_r>> J^r V\Sec
\\
        @VVV @VVV
\\
{S}(E,VJ^rQ) @>{S}(\varkappa_r)>> {S}(E,J^rVQ)
\end{CD}
\eEq
where the left and right arrows are the canonical injections
induced by (\ref{Eq10}) in combination with (\ref{Eq9}).

\section{The operator order on $J^r\F$}
\setcounter{equation}{0}
\setcounter{theorem}{0}

In \cite{CabKol95} there was discussed, in fact, the operator order
on an $F$-smooth morphism
\bEq\label{Eq16}
    {A}:\F\to \mc{F}(Y_1,Y)\,,
\eEq
where $Y\to M$ is another fibered manifold. We say that ${A}$ is
of the operator order $k$, if, for every
$\varphi,\psi\in C^\infty(Y_{1x},Y_{2x})$,
\beq
    j^k_y \varphi=j^k_y \psi \quad {\mathrm{implies}}\quad
        {A}(\varphi)(y) = {A}(\psi)(y)\,.
\eeq
Then ${A}$ determines the associated map
\bEq\label{Eq17}
    \mc{A} :\mc{F} J^k(Y_1,Y_2) \to Y\,, \quad
                \mc{A}(j^k_y \varphi) = A(\varphi)(y)\,,
\eEq
where
\beq
   \mc{F} J^k(Y_1,Y_2) =  \underset{x\in M}{\bigcup}
        J^k (Y_{1x},Y_{2x})
\eeq
is a classical manifold. By \cite{CabKol95}, $\mc{A}$ is a smooth map.

Let $x^i, x^p$ and $z^a$ be the local coordinates
on $Y_1$ and $Y_2$ from Section 2. Then the induced coordinates
on $\mc{F} J^k(Y_1,Y_2)$ are $z^a_\bet$,
$0\le \parallel\!\!\bet\!\!\parallel\le k$,
where $\bet$ is a multiindex of range $(m+1,\dots,m+n)$.
If $x^i, w^s$, $s=1,\dots,\dim Y-\dim M$, are local fiber
coordinates on $Y$, then the coordinate expression of $\mc{A}$ is
\beq
    w^s=f^s(x^i,x^p,z^a_\bet)\,.
\eeq

If we consider an $\mc{SB}$-morphism
\bEq\label{Eq18}
    A:J^r\F \to \mc{F}(Y_1,Y)\,,
\eEq
we have take into account that $\varphi,\psi\in J^r_x\F$ are
characterized by the associated maps
\beq
   \bar\varphi,\bar\psi:J^r_x Y_1\to J^r_x Y_2\,.
\eeq
So we have
$j^k\bar\varphi, j^k\bar\psi : J^r_x Y_1\to J^k(J^r_x Y_1,J^r_x Y_2)$.

\bDf
    We say that $A$ is of the operator order $k$, if
\beq
    j^k\bar\varphi|J^r_y Y_1=j^k\bar\psi|J^r_y Y_1\quad
{\mathrm{implies}}\quad A(\varphi)(y)=A(\psi)(y)\,.
\eeq
\eDf

To characterize the associated map of $A$ in this situation, we
introduce a new concept.

\bDf
    For a base preserving morphism $f:Y_1\to Y_2$, its
fiberwise $r$-jet prolongation $(\mc{F}j^r)f$ is defined by
\beq
    (\mc{F}j^r)f: Y_1\to \mc{F}J^r(Y_1,Y_2)\,,\quad
            (\mc{F}j^r)f(y) = j^r_y(f_x)\,,
\eeq
where $f_x:Y_{1x}\to Y_{2x}$ is the restricted and corestricted map,
$x=p_1(y)$. The $k$-jet $j^k_y (\mc{F}j^r)f$ is called
the {\em fiberwise $(k,r)$-jet} of $f$ at $y$.
\eDf

Let $\alp$ be a multiindex of the range $m$ and $\gam=(\alp,\bet)$.
Let $z^a= f^a(x^i,x^p)$ be the coordinate expression of $f$. Then
the coordinate expression of $(\mc{F}j^r)f$ is
\beq
    z^a_\bet= \der_\bet f^a\,,
        \quad 0\le\parallel\!\!\bet\!\!\parallel\le r\,.
\eeq
We write $\mc{F}J^{k,r}(Y_1,Y_2)= J^k(\mc{F}J^r(Y_1,Y_2)\to Y_1)$
for the space of all fiberwise $(k,r)$-jets of $Y_1$ to $Y_2$.
This is a classical manifold with the induced coordinates
\beq
    z^a_{\bet\gam}\,,\quad
        0\le\parallel\!\!\bet\!\!\parallel\le r\,,\,\,\,
        0\le\parallel\!\!\gam\!\!\parallel\le k\,.
\eeq
Clearly, we have
\bEq\label{Eq19}
    \mc{F} J^{k,0}(Y_1,Y_2) \simeq J^k(Y_1\ucar{M} Y_2\to Y_1)\,.
\eEq
Indeed, $(\mc{F}j^0)f=f$, which we identify with its graph
$Y_1\to Y_1\ucar{M} Y_2$, $y\mto (y,f(y))$.

\bPr\label{Pr2}
        If $A:J^r\F\to \mc{F}(Y_1,Y)$ is of operator order $k$, then
$A(j^r_x \widehat f)(y)$ depends on $j^k_y(\mc{F}j^r)f$ only.
\ePr

\bPf
    For $r=1$, the associated map $h: J^1_x Y_1\to J^1_x Y_2$
of an element of $J^1_x\F$ is
\bEq\label{Eq20}
    z^a= f^a(x^i_0, x^p),\,\,\, z^a_i=\der_i f^a(x^i_0,x^p) +
        \der_p f^a(x^i_0,x^p)\,x^p_i\,,
\eEq
where $x^p$ and $x^p_i$ are the variables on $J^1_x Y_1$,
$x=(x^i_0)\in M$.
Hence $j^k_y h|J^1_y Y$, $y=(x^i_0,x^p_0)$, depends on
\beq
    \der_\bet f^a(x^i_0,x^p_0),\,\,\,
\der_\bet\der_i f^a(x^i_0,x^p_0),\,\,\,
\der_\bet\der_p f^a(x^i_0,x^p_0),\,\,\,
0\le\parallel\!\!\bet\!\!\parallel\le k\,.
\eeq
These are the coordinates of $j^k_y(\mc{F}j^1)f$.
For $r>1$ we proceed by iteration using the facts $J^r$ is an $r$-th
order functor and the coordinate formula for $J^r f$ is of a specific
polynomial character in the induced jet coordinates.
\ePf

Hence $A$ determines the associated map
\beq
    \mc{A}: \mc{F}J^{k,r}(Y_1,Y_2)\to Y\,,
        \quad \mc{A}(j^k_y(\mc{F}j^r)f) = A(j^r_x\widehat f)(y)\,.
\eeq
Analogously to (\ref{Eq17}), $\mc{A}$ is a smooth map.

We remark that the concept of fiberwise $(k,r)$-jet can
be incorporated into the general framework of the concept of
$(r,s,q)$-jet of fibered manifold morphisms, \cite{KolMicSlo93}.
But this is somewhat sophisticated for our purposes, so that we prefer
our direct approach here.

\bRm\label{Rm1}
    It is interesting that a similar approach can be applied
to an arbitrary fiber product preserving bundle functor
$G$ on $\mc{FM}_m$. In \cite{CabJanKol04} we clarified that
$G$ can be extended to $\F$ as follows. If $G$ is of the base order $r$,
it can be identified with a triple $(A,H,t)$, where $A$ is a Weil
algebra, $H:G^r_m\to \Aut A$ is
a group homomorphism and $t:\mathbb{D}{}^r_m\to A$ is an equivariant
algebra homomorphism. In \cite{CabJanKol04} we defined
$G\F$ as the subset of the $F$-smooth associated bundle
$P^rM[T^A\F]$ of all equivariance classes $\{u,Z\}$,
$u\in P^rM,\,\, Z\in T^A\F$ satisfying $t_M(u)= T^Ap(Z)$.
\eRm

Analogously to the tangent case, $Z$ can be interpreted as a map
\beq
    \bar Z:T^A_XY_1 \to T^A_X Y_2\,,\quad X\in T^Ap(Z)\in T^AM\,.
\eeq
We know that $GY_i$, $i=1,2$, is the subset of $P^rM[T^AY_i]$ of
all $\{u,Z_i\}$ satisfying $t_M(u)= Tp_i(Z_i)$. Then
we construct a well defined inclusion
\beq
    G\F\subset \mc{F}(GY_1,GY_2)
\eeq
by transforming $\{u,Z\}\in G\F$ into the map
\beq
    \overline{\{u,Z\}}(\{u,Z_1\})=\{u,\bar Z(Z_1)\}\,,\quad \{u,Z_1\}
        \in GY_1\,.
\eeq
Thus, for every $G$ we can treat the operator order of an
$\mc{SB}$-morphism $G\F\to \mc{F}(Y_1,Y)$ similarly to the case
$G = J^r$.

\section{The formal exterior differential over $\F$}
\setcounter{equation}{0}
\setcounter{theorem}{0}

We recall that, given two other fibered manifolds $Y_3\to M$,
$Y_4\to M$, an $\mc{SB}$-morphism ${A}:\F\to \mc{F}(Y_3,Y_4)$
is called $J^k$-differentiable, if the rule
\beq
    (j^k_x\widehat s)\mto j^k_x({A}\com \widehat s),\quad \widehat s:M\to \F
\eeq
defines an $F$-smooth map
\beq
    J^k{A}:J^k\F\to J^k\mc{F}(Y_3,Y_4)\,.
\eeq

In general, consider three $F$-smooth bundles $S_1, S_2, S_3$ over $M$
and two surjective $\mc{SB}$-morphisms
$\pi_1:S_1\to S_3$, $\pi_2:S_2\to S_3$. We write
\beq
    S_1\ucar{S_3} S_2 = \{(u_1,u_2)\in S_1\car S_2,
        \pi_1(u_1)=\pi_2(u_2)\}\,\,.
\eeq
Clearly, this is also an $F$-smooth bundle over $M$.

Consider a $J^k$-differentiable morphism, $s\le r$,
\bEq\label{Eq21}
    A: J^r\F\ucar{J^s\F} VJ^s\F \to \mc{F}(Y_1,Y)\,.
\eEq
Using the exchange map $K_s$, see (\ref{Eq14}), we can define
\bEq\label{Eq22}
    \mc{J}^k_{\mathrm{hol}} A: J^{k+r}\F\ucar{J^{k+s}\F}
        VJ^{k+s}\F \to J^k\mc{F}(Y_1,Y)
\eEq
in the same way as in Section 1.

To introduce the formal exterior differential, we have to consider
${S}(Y_1,\bigwedge^l T^*Y_1)$ on the right hand side of (\ref{Eq21}).
So, let
\bEq\label{Eq23}
    A: J^r\F\ucar{J^s\F} VJ^s\F\to {S}(Y_1,\bigwedge^l T^*Y_1)
\eEq
be a $J^1$-differentiable morphism. Then we construct
$\mc{J}^1_{\mathrm{hol}} A$, use the inclusion
$$
i: J^1{S}(Y_1,\bigwedge^l T^* Y_1)\to {S}(Y_1, J^1\bigwedge^lT^* Y_1)
$$
on the right hand side and add
${S}(\del):{S}(Y_1,J^1\bigwedge^lT^* Y_1)\to{S}(Y_1,
\bigwedge^{l+1}T^* Y_1)$, where $\del$ is the formal version
of the exterior differential.

\bDf\label{Df4}
 The $F$-smooth morphism
\begin{align}\label{Eq24}
    \mathbb{D}A & \byd {S}(\del)\com i\com \mc{J}^1_{\mathrm{hol}} A :
\\
  & :
 J^{r+1}\F\ucar{J^{s+1}\F} VJ^{s+1}\F \to
        {S}(Y_1, \bigwedge^{l+1} T^* Y_1)\nonumber
\end{align}
will be called the {\em formal exterior differential} of (\ref{Eq23}).
\eDf

Clearly, the construction of $J^rY_1\ucar{J^sY_1} VJ^sY_1$ and of
the induced maps is a fiber product preserving bundle functor
on $\mc{FM}_m$. According to Remark \ref{Rm1}, we can introduce
the concept of operator order of the morphism (\ref{Eq21}).
However, we shall not go into details in this paper.

\section{The restriction of $\mathbb D$ to $\Sec$}
\setcounter{equation}{0}
\setcounter{theorem}{0}

Now we consider $\Sec$ in the role of $\F$. Let $P\to E\to M$
be another 2-fibered manifold and
\bEq\label{Eq25}
{A}:J^r\Sec\ucar{J^s\Sec}VJ^s\Sec\to {S}(E,P)
\eEq
be a $J^k$-differentiable morphism. Then (\ref{Eq22}) restricts
to a morphism
\bEq\label{Eq26}
    \mc{J}^k_{\mathrm{hol}} A: J^{k+s}\Sec \ucar{J^{k+s}\Sec}
        VJ^{k+s}\Sec \to J^k{S}(E,P)\,.
\eEq

In the case $k=1$ and $P=\bigwedge^l T^* E$, (\ref{Eq24}) yields
a morphism
\begin{align}\label{Eq27}
    \mathbb{D}A & :
 J^{r+1}\Sec\ucar{J^{s+1}\Sec} VJ^{s+1}\Sec \to
        {S}(E, \bigwedge^{l+1} T^* E)\,.
\end{align}

An important fact is that (\ref{Eq27}) and the finite dimensional
formal exterior differential over $E$ are related as follows. From
now on we always consider $Q$ as a fibered manifold over $E$. Let
\bEq\label{Eq28}
    B: J^{r}Q \ucar{J^{s} Q}
        VJ^{s}Q \to P\,
\eEq
be a smooth $E$-morphism. On one hand, we construct
\bEq\label{Eq29}
    {S}(B): {S}(E,J^{r}Q) \ucar{{S}(E,J^{s} Q)}
        {S}(E,VJ^{s}Q) \to {S}(E,P)\,.
\eEq
The injection $J^r\Sec\to {S}(E,J^rQ)$ induces,
including holonomization, an injection
\begin{align}\label{Eq30}
    I : J^k{S}(E,J^{r}Q) & \ucar{J^k{S}(E,J^{s} Q)}
        J^k {S}(E,VJ^{s}Q)
        \to
\\
&\to
        {S}(E,J^{k+r} Q)\ucar{{S}(E,J^{k+s} Q)} {S}(E,VJ^{k+s} Q)\,.
        \nonumber
\end{align}
On the other hand, we can construct
\bEq\label{Eq31}
    \mc{J}^k_{\mathrm{hol}}B: J^{k+r}Q \ucar{J^{k+s} Q}
        VJ^{k+s}Q \to J^k P\,.
\eEq
So we have a diagram
\bEq\label{Eq32}
\begin{CD}
    J^k{S}(E,J^rQ)\ucar{J^k{S}(E,J^sQ)} J^k{S}(E,VJ^s Q)
        @>J^k{S}(B)>> J^k {S}(E,P)
\\
        @VIVV @ViVV
\\
{S}(E,J^{k+r}Q)\ucar{{S}(E,J^{k+s} Q)} {S}(E,VJ^{k+s} Q)
        @>{S}(\mc{J}^k_{\mathrm{hol}}B)>> {S}(E,J^k P)
\end{CD}
\eEq
Then the proofs of (\ref{Eq13}) and (\ref{Eq15}) imply

\bLm\label{Lm2}
(\ref{Eq32}) is a commutative diagram. \hfill {\qedsymbol}
\eLm

In the case of $B: J^{r}Q \ucar{J^{s} Q}
        VJ^{s}Q \to \bigwedge^l T^* E$,
${S}(B)$ induces
\bEq\label{Eq33}
    \mathbb{D}{S}(B) : J^1{S}(E,J^{r}Q) \ucar{J^1{S}(E,J^{s} Q)}
        J^1 {S}(E,VJ^{s}Q)
        \to
        {S}(E,\bigwedge^{l+1} T^* E)\,.
\eEq
On the other hand, we have
${D}B:J^{r+1} Q\car_{J^{s+1}Q} VJ^{s+1} Q\to \bigwedge^{l+1} T^* E$.
Then Lemma \ref{Lm2} implies

\bPr\label{Pr3}
  We have $\mathbb{D}({S}(B)) = {S}({D}B)\com I$.
\hfill {\qedsymbol}
\ePr

\section{The Euler-Lagrange morphism}
\setcounter{equation}{0}
\setcounter{theorem}{0}

We first recall a suitable construction of the Euler-Lagrange morphism
of a first order Lagrangian on a fibered manifold $Y\to M$,
\cite{Kol84, KolVit03}. We shall discuss a slightly more general case
of a morphism
\bEq\label{Eq34}
    \lam : J^1 Y
        \to
       \bigwedge^{l} T^* M\,.
\eEq
If $l=m=\dim M$, we obtain a classical first order Lagrangian on $Y$.

The vertical differential of $\lam$ is a map
\bEq\label{Eq35}
    d_V\lam : J^1 Y
        \to
       V^*J^1Y\ten \bigwedge^{l} T^* M\,.
\eEq
The well-known exact sequence
\beq
    0\to VY\ten T^* M\to VJ^1Y \to VY \to 0
\eeq
induces the dual map $V^*J^1 Y\to V^* Y\ten TM$. If we add
the classical tensor contraction
$\text{\aa \char'171} : TM\ten \bigwedge^lT^* M\to \bigwedge^{l-1} T^* M$,
we obtain the composed map
\bEq\label{Eq36}
    \rho_Y : V^*J^1 Y \ten \bigwedge^l T^* M
        \to
       V^*Y\ten \bigwedge^{l-1} T^* M\,.
\eEq
Hence $\rho_Y\com d_V\lam$ can be interpreted as a morphism
\bEq\label{Eq37}
    B(\lam) = \rho_Y\com d_V\lam : J^1 Y \ucar{Y} VY \to
         \bigwedge^{l-1} T^* M\,.
\eEq
Then
\beq
{D}B(\lam) : J^2 Y \ucar{J^1 Y} VJ^1 Y \to
         \bigwedge^{l} T^* M\,.
\eeq
In coordinates, if
\beq
    \lam \equiv L_{i_1\dots i_l}(x^i,x^p,x^p_i)\,
        dx^{i_1}\wed\dots \wed dx^{i_l}\,,
\eeq
then
\beq
    d_V\lam \equiv \frac{\der L_{i_1\dots i_l}}{\der x^p}\,
        dx^p\ten dx^{i_1}\wed\dots \wed dx^{i_l}  +
        \frac{\der L_{i_1\dots i_l}}{\der x^p_i}\,
        dx^p_i\ten dx^{i_1}\wed\dots \wed dx^{i_l}  \,.
\eeq
Hence
\beq
    B(\lam)\equiv  \frac{\der L_{i_1\dots i_l}}{\der x^p_i}\, dx^p\ten
        \frac{\der}{\der x^i} \, \text{\aa \char'171} \,
        dx^{i_1}\wed\dots \wed dx^{i_l}  \,.
\eeq
Using (\ref{Eq1.9}), we obtain
\beq
{D}B(\lam)\equiv D_i \frac{\der L_{i_1\dots i_l}}{\der x^p_i}\,
        dx^p\ten dx^{i_1}\wed\dots \wed dx^{i_l} +
        \frac{\der L_{i_1\dots i_l}}{\der x^p_i}\, dx^p_i\ten
        dx^{i_1}\wed\dots \wed dx^{i_l}\,,
\eeq
where $D_i$ denotes the standard formal partial derivative
with respect to $x^i$, \cite{Kol84}. This implies that the difference
\bEq\label{Eq38}
    \mc{E}(\lam) \byd (d_V\lam)\com \pi^2_1 - {D} B(\lam)\,,
\eEq
where $\pi^2_1:J^2 Y\to J^1 Y$ is the jet map, is projectable to $VY$.
Hence it can be interpreted as a morphism
\bEq\label{Eq38}
    \mc{E}(\lam): J^2 Y\to V^* Y\ten \bigwedge^l T^* M\,.
\eEq
Its coordinate expression implies that for $m=l$ we obtain
the Euler-Lagrange morphism of $\lam$. A more geometric explanation
of this fact can be found in \cite{Kol84}.

Consider now a smooth $E$-morphism
\bEq\label{Eq39}
    L : J^1 Q\to  \bigwedge^l T^* E\,,
\eEq
which is a first order Lagrangian on $Q\to E$ in the case $l=m+n$.
We can interpret $L$ as an $\mc{SB}$-morphism
\bEq\label{Eq40}
    \widehat{L} : {S}(E,J^1 Q)\to {S}(E, \bigwedge^l T^* E)\,
\eEq
or as a section, denoted by the same symbol,
\bEq\label{Eq41}
    \widehat{L} : M\to \mc{F}(J^1 Q, \bigwedge^l T^* E)\,.
\eEq
Under the functional approach, the vertical differential is a map
\bEq\label{Eq42}
    \widehat{d}_V : \mc{F}(J^1 Q, \bigwedge^l T^* E)
        \to {S}(J^1 Q, V^*J^1Q\ten \bigwedge^l T^* E)\,.
\eEq
Then we take into account
\beq
\rho_Q : V^*J^1 Q \ten \bigwedge^l T^* E
        \to
       V^*Q\ten \bigwedge^{l-1} T^* E\,.
\eeq
This induces, fiberwise, a map
\beq
    \mc{F}(\id_{J^1Q},\rho_Q):
        \mc{F}(J^1Q, V^*J^1Q\ten \bigwedge^lT^* E)
        \to \mc{F}(J^1 Q, V^*Q\ten \bigwedge^{l-1} T^* E)\,.
\eeq
Then $\widehat{B}(\widehat{L})\byd \mc{F}(\id_{J^1Q},\rho_Q)\com
\widehat{d}_V\com \widehat{L}$ can be viewed as a morphism
\beq
\widehat{B}(\widehat{L}):{S}(E,J^1 Q\ucar{Q} VQ)
        \to {S}(E,\bigwedge^{l-1} T^* E)
\eeq
and we can construct
\beq
\mathbb{D}(\widehat{B}(\widehat{L})):J^1 {S}(E,J^1 Q)
        \ucar{J^1{S}(E,Q)} J^1{S}(E,VJ^1 Q)
        \to {S}(E,\bigwedge^{l} T^* E)\,.
\eeq
Since $\widehat{B}(\widehat{L})= {S}(B(L))$, Proposition \ref{Pr3}
yields $\mathbb{D}(\widehat{B}(\widehat{L}))={S}(D(B(L)))\com I$.
This implies that $\mathbb{D}(\widehat{B}(\widehat{L}))$ can be
viewed as an $F$-smooth section
\beq
    M\to \mc{F}(J^2 Q,V^*J^1 Q\ten \bigwedge^l T^*E)
\eeq
and
$\widehat{d}_V\com \widehat{L}\com \pi^2_1
-\mathbb{D}(\widehat{B}(\widehat{L}))$
corestricts to an $F$-smooth section
\bEq\label{Eq43}
\widehat{\mc{E}}(\widehat{L})=
\widehat{d}_V\com \widehat{L}\com \pi^2_1
-\mathbb{D}(\widehat{B}(\widehat{L})):
        M\to \mc{F}(J^2 Q, V^* Q\ten \bigwedge^lT^* E)\,.
\eEq
This can be viewed as an $\mc{SB}$-morphism, denoted by the same symbol,
\bEq\label{Eq44}
\widehat{\mc{E}}(\widehat{L}):
        {S}(E,J^2 Q)\to {S}(E, V^* Q\ten \bigwedge^lT^* E)\,.
\eEq
By construction, $\widehat{\mc{E}}(\widehat{L})={S}(\mc{E}(L))$.

Thus, in the case $l=m+n$ our construction represents a functional
approach to the Euler-Lagrange morphism of a first order Lagrangian
on $Q\to E$.

\bRm\label{Rm2}
Given two fibered manifolds $Y_1\to M$ and $Y_2\to M$,
one can study variational calculus for the base preserving morphisms
$Y_1\to Y_2$. Since these morphisms are identified with the sections
of the fibered manifold $Y_1\car_{M} Y_2\to Y_1$, from the abstract
point of view the morphism problem reduces to the variational calculus
for the sections of the latter bundle. However, the geometry
of the morphism problem should be more rich. It seems to be
reasonable to discuss this subject in more details elsewhere.
\eRm



\end{document}